\documentclass{article}
\usepackage{amssymb,latexsym}
\usepackage{amsmath}
\usepackage{graphics}
\usepackage{graphicx}
\newtheorem{theorem}{Theorem}
\newtheorem{definition}{Definition}

\newtheorem{remark}{Remark}

\numberwithin{equation}{section} \numberwithin{theorem}{section}
\numberwithin{lemma}{section} \numberwithin{remark}{section}
\numberwithin{notation}{section}
\numberwithin{definition}{section}

\begin{document}
\title{Generalized Mannheim Curves in Minkowski space-time $E_1^4$}
\author{Soley Ersoy$^{a}$ , Murat Tosun$^{a}$ , Hiroo Matsuda$^{b}$ \\
{\small sersoy@sakarya.edu.tr , tosun@sakarya.edu.tr ,
matsuda@kanazawa-med.ac.jp }\\
 {\small {$^{a}$ Department of
Mathematics, Sakarya University,
Sakarya, TURKEY} }\\
{\small {$^{b}$ Department Mathematics, Kanazawa Medical
University, Uchinada, Ishikawa, 920-02, JAPAN} }}

\date{}

\maketitle

\maketitle

\begin{abstract}
In this paper, the definition of generalized spacelike Mannheim
curve in Minkowski space-time $E_1^4$ is given. The necessary and
sufficient conditions for the generalized spacelike Mannheim curve
are obtained. Also, some characterizations of Mannheim curve are
given.

\textbf{Mathematics Subject Classification (2010):} 53B30, 53A35,
53A04.

\textbf{Keywords}: Mannheim curve, Minkowski space-time\\
\end{abstract}

\section{Introduction}\label{S:intro}

The curves are a fundamental structure of differential geometry.
An increasing interest of the theory of curves makes a development
of special curves to be examined. A way to classification and
characterization of curves is the relationship between the Frenet
vectors of the curves. For example, Saint Venant proposed the
question whether upon the surface generated by the principal
normal of a curve, a second curve can exist which has for its
principal normal of the given curve in 1845. This question was
answered by Bertrand in 1850. He showed that a necessary and
sufficient condition for the existence of such a second curve is
that a linear relationship with constant coefficients exists
between the first and second curvatures of the given original
curve. The pairs of curves of this kind have been called Bertrand
partner curves or more commonly Bertrand curves \cite{Car},
\cite{Kuh}, \cite{Str}. There are many works related with Bertrand
curves in the Euclidean space and Minkowski space,
\cite{Nad}--\cite{Bal2}. Also, generalized Bertrand curves in
Euclidean 4- space are defined and characterized in \cite{Mat1}.
Another kind of associated curve have been called Mannheim curve
and Mannheim partner curve. The notion of Mannheim curves was
discovered by A. Mannheim in 1878. These curves in Euclidean
3-space are characterized in terms of the curvature and torsion as
follows: A space curve is a Mannheim curve if and only if its
curvature $\kappa $ and torsion $\tau$ satisfy the relation
$$
\kappa \left( s \right) = \alpha \left( {{\kappa ^2}\left( s
\right) + {\tau ^2}\left( s \right)} \right)
$$
for some constant $\alpha$. The articles concerning Mannheim
curves are rather few. In \cite{Blum}, a remarkable class of
Mannheim curves is studied. General Mannheim curves in the
Euclidean 3-space are obtained in \cite{Tig}. Mannheim partner
curves in Euclidean 3-space and Minkowski 3-space are studied and
the necessary and sufficient conditions for the Mannheim partner
curves are obtained in \cite{Liu}, \cite{Orb}. Recently, Mannheim
curves are generalized and some characterizations and examples of
generalized Mannheim curves in Euclidean 4-space ${E^4}$ are given
by \cite{Mat2}.

\noindent In this paper, we study the generalized spacelike
Mannheim partner curves in $4-$dimensional Minkowski space-time.
We will give the necessary and sufficient conditions for the
generalized spacelike Mannheim partner curves.

\section{Preliminaries}\label{S:intro}

The basic concepts of the theory of curves in Minkowski space-time
${E^4}$ are briefly presented in this section. A more complete
elementary treatment can be found in \cite{Onei}. Minkowski
space-time ${E_1^4}$ is an Euclidean space provided with the
standard flat metric given by
$$
\left\langle {\,\,,\,\,} \right\rangle  =  - dx_1^2 + dx_2^2 +
dx_3^2 + dx_4^2
$$
where $\left( {{x_1},\,{x_2},\,{x_3},\,{x_4}} \right)$ is a
rectangular coordinate system in ${E^4}$.

\noindent Since $\left\langle {\;,\;} \right\rangle $ is an
indefinite metric, recall that a vector ${\bf{v}} \in E_1^4$ can
have one of the three causal characters; it can be spacelike if
$\left\langle {{\bf{v}},{\bf{v}}} \right\rangle  > 0$ or ${\bf{v}}
= {\bf{0}}$, timelike if $\left\langle {{\bf{v}},{\bf{v}}}
\right\rangle  < 0$ and null (lightlike) if $\left\langle
{{\bf{v}},{\bf{v}}} \right\rangle  = 0$ and ${\bf{v}} \ne
{\bf{0}}$ . Similarly, an arbitrary curve ${\bf{c}} =
{\bf{c}}\left( s \right)$  in ${E^4}$ can locally be spacelike,
timelike or null (lightlike) if all of its velocity vectors
${\bf{c'}}\left( s \right)$ are, respectively, spacelike, timelike
or null. The norm of ${\bf{v}} \in E_1^4$ is given by $\left\|
{\bf{v}} \right\| = \sqrt {\left| {\left\langle
{{\bf{v}},{\bf{v}}} \right\rangle } \right|} $. If $\left\|
{{\bf{c'}}\left( s \right)} \right\| = \sqrt {\left| {\left\langle
{{\bf{c'}}\left( s \right),{\bf{c'}}\left( s \right)}
\right\rangle } \right|}  \ne 0$ for all $s \in L$, then $C$ is a
regular curve in $E_1^4$. A spacelike (timelike) regular curve $C$
is parameterized by arc-length parameter $s$ which is given by
${\bf{c}}:L \to E_1^4$, then the tangent vector ${\bf{c'}}\left( s
\right)$ along $C$ has unit length, that is, \[\left\langle
{{\bf{c}}\left( s \right),{\bf{c}}\left( s \right)} \right\rangle
= 1\,,\,\quad \left( {\left\langle {{\bf{c}}\left( s
\right),{\bf{c}}\left( s \right)} \right\rangle  =  - 1} \right)\]
for all $s \in L$ .

\noindent Hereafter, curves are considered  spacelike and regular
${C^\infty}$ curves in $E_1^4$. Let ${{\bf{e}}_1}\left( s \right)
= {\bf{c'}}\left( s \right)$ for all $s \in L$, then the vector
field ${{\bf{e}}_1}\left( s \right)$ is spacelike and it is called
spacelike unit tangent vector field on $C$.

\noindent The spacelike curve $C$ is called special spacelike
Frenet curve if there exist three smooth functions ${k_1}$,
${k_2}$, ${k_3}$ on $C$ and smooth non-null frame field
$\left\{{{{\bf{e}}_1},\,{{\bf{e}}_2},\,{{\bf{e}}_3},\,{{\bf{e}}_4}}
\right\}$ along the curve $C$. Also, the functions
${k_1},\,{k_2}$, and ${k_3}$ are called the first, the second, and
the third curvature function on $C$, respectively. For the
${C^\infty }$ special spacelike Frenet curve $C$, the following
Frenet formula is hold

\[\left[ \begin{array}{l}
 {{{\bf{e'}}}_1} \\
 {{{\bf{e'}}}_2} \\
 {{{\bf{e'}}}_3} \\
 {{{\bf{e'}}}_4} \\
 \end{array} \right] = \,\left[ \begin{array}{l}
 \,\,\,0\,\,\,\,\,\,\,\,\,\,{k_1}\,\,\,\,\,\,\,\,\,0\,\,\,\,\,\,\,\,\,0 \\
 {\mu _1}{k_1}\,\,\,\,\,\,\,0\,\,\,\,\,\,\,\,\,\,{k_2}\,\,\,\,\,\,\,0 \\
 \,\,\,0\,\,\,\,\,\,\,{\mu _2}{k_2}\,\,\,\,\,\,\,0\,\,\,\,\,\,\,\,\,{k_3} \\
 \,\,\,0\,\,\,\,\,\,\,\,\,\,\,0\,\,\,\,\,\,\,\,{\mu _3}{k_3}\,\,\,\,0 \\
 \end{array} \right]\left[ \begin{array}{l}
 {{\bf{e}}_1} \\
 {{\bf{e}}_{\rm{2}}} \\
 {{\bf{e}}_{\rm{3}}} \\
 {{\bf{e}}_4} \\
 \end{array} \right]\]
where ${\mu _i} =  \mp 1,\,\,1 \le i \le 3$, \cite{Onei}.\\
Due to characters of Frenet vectors of the spacelike curve $C$,
${\mu _i}\,\,\left( {1 \le i \le 3} \right)$ are defined as in the
following three subcases; \noindent \textbf{Case 1:} If
${{\bf{e}}_4}$ is timelike, then ${\mu _i},\,\,1 \le i \le 3$ are
$$
{\mu _1} = {\mu _2} =  - 1\,\,,\,\,\,{\mu _3} = 1
$$
where ${{\bf{e}}_1},\,{{\bf{e}}_2},\,{{\bf{e}}_3}$ and
${{\bf{e}}_4}$ are mutually orthogonal vector fields satisfying
equations
$$
\left\langle {{{\bf{e}}_1}\,,\,{{\bf{e}}_1}} \right\rangle  =
\left\langle {{{\bf{e}}_2}\,,\,{{\bf{e}}_2}} \right\rangle  =
\left\langle {{{\bf{e}}_3}\,,\,{{\bf{e}}_3}} \right\rangle  =
1\,\,,\,\,\left\langle {{{\bf{e}}_4}\,,\,{{\bf{e}}_4}}
\right\rangle  =  - 1.
$$

\noindent \textbf{Case 2:} If ${{\bf{e}}_3}$ is timelike, then
${\mu _i},\,\,1 \le i \le 3$ are
$$
{\mu _1} =  - 1\,\,,\,\,{\mu _2} = {\mu _3} = 1
$$
where ${{\bf{e}}_1},\,{{\bf{e}}_2},\,{{\bf{e}}_3}$ and
${{\bf{e}}_4}$ are mutually orthogonal vector fields satisfying
equations
$$
\left\langle {{{\bf{e}}_1}\,,\,{{\bf{e}}_1}} \right\rangle  =
\left\langle {{{\bf{e}}_2}\,,\,{{\bf{e}}_2}} \right\rangle  =
\left\langle {{{\bf{e}}_4}\,,\,{{\bf{e}}_4}} \right\rangle  =
1\,\,,\,\,\left\langle {{{\bf{e}}_3}\,,\,{{\bf{e}}_3}}
\right\rangle  =  - 1.
$$

\noindent \textbf{Case 3:} If ${{\bf{e}}_2}$ is timelike, then
${\mu _i},\,\,1 \le i \le 3$ are
$$
{\mu _1} = {\mu _2} = 1\,\,,\,\,{\mu _3} =  - 1
$$
where ${{\bf{e}}_1},\,{{\bf{e}}_2},\,{{\bf{e}}_3}$ and
${{\bf{e}}_4}$ are mutually orthogonal vector fields satisfying
equations
$$
\left\langle {{{\bf{e}}_1}\,,\,{{\bf{e}}_1}} \right\rangle  =
\left\langle {{{\bf{e}}_3}\,,\,{{\bf{e}}_3}} \right\rangle  =
\left\langle {{{\bf{e}}_4}\,,\,{{\bf{e}}_4}} \right\rangle  =
1\,\,,\,\,\,\left\langle {{{\bf{e}}_2}\,,\,{{\bf{e}}_2}}
\right\rangle  =  - 1.
$$

\noindent For $s \in L$, the non-null frame field $\left\{
{{{\bf{e}}_1},\,{{\bf{e}}_2},\,{{\bf{e}}_3},\,{{\bf{e}}_4}}
\right\}$ and curvature functions ${k_1}$ and ${k_2}$ are
determined as follows

\[\begin{array}{l}
 {1^{st}}\,\,\,\,\,\,{\rm{step    }}\,\,\,\,\,\,{{\bf{e}}_1}\left( s \right) = {\bf{c'}}\left( s \right) \\
 {2^{nd}}\,\,\,\,{\rm{step   }}\,\,\,\,\,\,\,{k_1}\left( s \right) = \left\| {{{{\bf{e'}}}_1}\left( s \right)} \right\| > 0 \\
 \,\,\,\,\,\,\,\,\,\,\,\,\,\,\,\,\,\,\,\,\,\,\,\,\,\,\,\,\,\,\,{{\bf{e}}_2}\left( s \right) = \frac{1}{{{k_1}\left( s \right)}}{{{\bf{e'}}}_1}\left( s \right) \\
 {3^{rd}}\,\,\,\,{\rm{step}}\,\,\,\,\,\,\,\,{k_2}\left( s \right) = \left\| {{{{\bf{e'}}}_2}\left( s \right) - {\mu _1}{k_1}\left( s \right){{\bf{e}}_1}\left( s \right)} \right\| > 0 \\
 \,\,\,\,\,\,\,\,\,\,\,\,\,\,\,\,\,\,\,\,\,\,\,\,\,\,\,\,\,\,\,\,{{\bf{e}}_3}\left( s \right) = \frac{1}{{{k_2}\left( s \right)}}\left( {{{{\bf{e'}}}_2}\left( s \right) - {\mu _1}{k_1}\left( s \right){{\bf{e}}_1}\left( s \right)} \right) \\
 {4^{th}}\,\,\,{\rm{step}}\,\,\,\,\,\,\,\,\,\,{e_4}\left( s \right) = \varepsilon \frac{1}{{\left\| {{{{\bf{e'}}}_3}\left( s \right) - {\mu _2}{k_2}\left( s \right){{\bf{e}}_2}\left( s \right)} \right\|}}\left( {{{{\bf{e'}}}_3}\left( s \right) - {\mu _2}{k_2}\left( s \right){{\bf{e}}_2}\left( s \right)} \right) \\
 \end{array}\]
where $\varepsilon $ is taken $ - 1$ or $ +1$ to make $ + 1$ the
determinant of $\left\{
{{{\bf{e}}_1},\,{{\bf{e}}_2},\,{{\bf{e}}_3},\,{{\bf{e}}_4}}
\right\}$, that is, the non-null orthonormal frame field is of
positive orientation. The function ${k_3}$ is determined by
\[{k_3}\left( s \right) = \left\langle {{{{\bf{e'}}}_3}\left( s \right)\,,\,{{\bf{e}}_4}\left( s \right)} \right\rangle  \ne 0.\]
So the function ${k_3}$ never vanishes.

\noindent In order to make sure that the spacelike curve $C$ is a
special spacelike Frenet curve, above steps must be checked, from
${1^{st}}$ step to ${4^{th}}$ step, for $s \in L$.

\noindent At each point of spacelike curve $C$, a line ${\ell _1}$
in the direction of ${{\bf{e}}_2}$ is called the first normal
line, a line ${\ell _2}$ in the direction of ${{\bf{e}}_3}$ is
called the second normal line and a line ${\ell _3}$ in the
direction of ${{\bf{e}}_4}$ is called the third normal line.

\noindent Note that, according to three different case of
spacelike curve $C$, ${\ell _3},\,{\ell _2}$ and ${\ell _1}$ can
be timelike, respectively, which are called second binormal, first
binormal and principal normal line at each point of the spacelike
curve $C$.

\section{Generalized spacelike Mannheim curves in  $E_1^4$}\label{S:intro}

In ${E^4}$ the Bertrand curves and Mannheim curves  are
generalized by \cite{Mat1} and \cite{Mat2}, respectively. In these
regards, we have investigate generalization of spacelike Mannheim
curves  Minkowski space in $E_1^4$.
\begin{definition}
A special spacelike curve $C$ in $E_1^4$ is a generalized
spacelike Mannheim curve if there exists a special spacelike
Frenet curve ${C^*}$ in $E_1^4$ such that the first normal line at
each of $C$ is included in the plane generated by the second
normal line and the third normal line of ${C^*}$ at the
corresponding point under $\phi $. Here $\phi $ is a bijection
from $C$ to ${C^*}$. The curve  ${C^*}$ is called the generalized
spacelike Mannheim mate curve of $C$.
\end{definition}

\noindent By the definition, a generalized Mannheim mate curve
${C^*}$ is given by

\begin{equation}\label{3.1}
\begin{array}{l}
{{\bf{c}}^*}\left( s \right) = {\bf{c}}\left( s \right) + \alpha
\left( s \right){{\bf{e}}_2}\left( s \right),\,\,s \in L
 \end{array}
\end{equation}
where $\alpha $ is a smooth function on $L$. Generally, the
parameter $s$ isn't an arc-length of ${C^*}$. Let ${s^*}$ be the
arc-length of ${C^*}$ defined by

$${s^*} = \int\limits_0^s {\left\| {\frac{{d{{\bf{c}}^*}\left( s
\right)}}{{ds}}} \right\|ds.}
$$
If a smooth function $f:L \to L$ is given by $f\left( s \right) =
{s^*}$, then

$$
\begin{array}{l}
 \,\,\,\,\frac{{d{{\bf{c}}^*}\left( s \right)}}{{ds}} = \,\,\,{{\bf{e}}_1}\left( s \right) + \alpha '\left( s \right){{\bf{e}}_2}\left( s \right) + \alpha \left( s \right){\mu _1}{k_1}\left( s \right){{\bf{e}}_1}\left( s \right) + \alpha \left( s \right){k_2}\left( s \right){{\bf{e}}_3}\left( s \right) \\
 \,\,\,\,\,\,\,\,\,\,\,\,\,\,\,\,\,\,\;\, = \,\,\,\left( {1 + {\mu _1}\alpha \left( s \right){k_1}\left( s \right)} \right){{\bf{e}}_1}\left( s \right) + \alpha '\left( s \right){{\bf{e}}_2}\left( s \right) + \alpha \left( s \right){k_2}\left( s \right){{\bf{e}}_3}\left( s \right). \\
 \end{array}
$$
for $\forall s \in L$. Thus, we have
$$
\begin{array}{l}
f'\left( s \right) = \frac{{d{s^*}}}{{ds}} = \left\|
{\frac{{d{{\bf{c}}^{\bf{*}}}\left( s \right)}}{{ds}}} \right\| =
\sqrt {\left| {{{\left( {1 + {\mu _1}\alpha \left( s
\right){k_1}\left( s \right)} \right)}^2} + \varepsilon_2 {{\left(
{\alpha '\left( s \right)} \right)}^2} + \varepsilon_3 {{\left(
{\alpha \left( s \right){k_2}\left( s \right)} \right)}^2}}
\right|}
 \end{array}
$$
where $\varepsilon_i  = \left\{ \begin{array}{l}
  - 1\,\,,\,\,\,{{\bf{e}}_i}\,\,{\rm{is }}\,{\rm{timelike}} \\
 \,\,\,\,1\,\,,\,\,\,\,{{\bf{e}}_i}\,{\rm{is }}\,{\rm{spacelike}} \\
 \end{array} \right.$,  for $2 \le i \le 4.$
\\This means that, in the Case 1, ${{\bf{e}}_4}$ is timelike and

$$
\begin{array}{l}
f'\left( s \right) = \sqrt {\left| {{{\left( {1 - \alpha \left( s
\right){k_1}\left( s \right)} \right)}^2} + {{\left( {\alpha
'\left( s \right)} \right)}^2} + {{\left( {\alpha \left( s
\right){k_2}\left( s \right)} \right)}^2}} \right|}
 \end{array}
$$
or in the Case 2, ${{\bf{e}}_3}$ is timelike and
$$
\begin{array}{l}
f'\left( s \right) = \sqrt {\left| {{{\left( {1 - \alpha \left( s
\right){k_1}\left( s \right)} \right)}^2} + {{\left( {\alpha
'\left( s \right)} \right)}^2} - {{\left( {\alpha \left( s
\right){k_2}\left( s \right)} \right)}^2}} \right|}
 \end{array}
$$
or  in the Case 3, ${{\bf{e}}_2}$ is timelike and
$$
\begin{array}{l}
f'\left( s \right) = \sqrt {\left| {{{\left( {1 + \alpha \left( s
\right){k_1}\left( s \right)} \right)}^2} - {{\left( {\alpha
'\left( s \right)} \right)}^2} + {{\left( {\alpha \left( s
\right){k_2}\left( s \right)} \right)}^2}} \right|}.
 \end{array}
$$
The spacelike curve ${C^*}$ with arc-length parameter ${s^*}$ is
$$
\begin{array}{l}
 {{\bf{c}}^*}:\,{L^*} \to E_1^4 \\
 \,\,\,\,\,\,\,\,\,\,\,{s^*}\,\, \to \,{{\bf{c}}^*}\left( {{s^*}} \right). \\
 \end{array}
$$
For a bijection $\phi :\,C \to {C^*}$ defined by $\phi
\left({{\bf{c}}\left( s \right)} \right) = {{\bf{c}}^*}\left(
{f\left( s \right)} \right),$ the reparametrization of ${C^*}$ is
$$
\begin{array}{l}
{{\bf{c}}^*}\left( {f\left( s \right)} \right) = {\bf{c}}\left( s
\right) + \alpha \left( s \right){{\bf{e}}_2}\left( s \right)
 \end{array}
$$
where $\alpha $ is a smooth function on $L$.

\begin{theorem}\label{T:3.1}
If a special spacelike Frenet curve $C$ in $E_1^4$ is a
generalized spacelike Mannheim curve, then the first curvature
function ${k_1}$ and the second curvature function
 ${k_2}$ of $C$ satisfy the equality
\begin{equation}\label{3.2}
\begin{array}{l}
{k_1}\left( s \right) =  - \alpha \left( {{\mu _1}k_1^2\left( s
\right) + {\mu _2}k_2^2\left( s \right)} \right)\,\,,\,\,s \in L
 \end{array}
\end{equation}
where $\alpha $ is a constant number and ${\mu _1} = {\mu _2} =  -
1$ when ${{\bf{e}}_4}$ is timelike or ${\mu _1} =  - 1\,,\,\,{\mu
_2} = 1$ when ${{\bf{e}}_3}$ is timelike or ${\mu _1} = {\mu _2} =
1$ when ${{\bf{e}}_2}$ is timelike.
\end{theorem}

\noindent \textbf{Proof.} Let $C$ be a generalized spacelike
Mannheim curve and ${C^*}$ be the generalized spacelike Mannheim
mate curve of $C$ with the diagram;
$$
\begin{array}{l}
 \,\,\,\,\,\,\,\,\,\,\,\,\mathop {\bf{c}}\limits_{ \cdot \,\, \cdot } \,\,\,\,\,\,\,\,\,\,\,\,\,\,\,\,\,\,\,{\mathop {\bf{c}}\limits_{ \cdot \,\, \cdot } ^*}\, \\
 f:\,\,\,\,L\,\,\,\,\,\, \to \,\,\,\,\,{L^*}\, \\
 \,\,\,\,\,\,\,\,\,\,\,\,\,\, \downarrow \,\,\,\,\,\,\,\,\,\,\,\,\,\,\,\,\,\,\,\, \downarrow  \\
 \phi \,:\,\,\,E_1^4\,\,\, \to \,\,\,E_1^4. \\
 \end{array}
$$
A smooth function $f$ is defined by $f\left( s \right) = \int
{\left\| {\frac{{d{{\bf{c}}^*}\left( s \right)}}{{ds}}} \right\|}
ds = {s^*}$ and ${s^*}$ is the arc-length parameter of ${C^*}$.
Also $\phi $ is a bijection which is defined by $\phi \left(
{{\bf{c}}\left( s \right)} \right) = {{\bf{c}}^*}\left( {f\left( s
\right)} \right).$ Thus, the spacelike curve ${C^*}$ is
reparametrized by
\begin{equation}\label{3.3}
\begin{array}{l}
{{\bf{c}}^*}\left( {f\left( s \right)} \right) = {\bf{c}}\left( s
\right) + \alpha \left( s \right){{\bf{e}}_2}\left( s \right)
 \end{array}
\end{equation}
where $\alpha $ is a smooth function. By differentiating both
sides of (3.3) with respect to $s$
\begin{equation}\label{3.4}
\begin{array}{l}
 f'\left( s \right){\bf{e}}_1^*\left( {f\left( s \right)} \right) = \left( {1 + {\mu _1}\alpha \left( s \right){k_1}\left( s \right)} \right){{\bf{e}}_1} + \alpha '\left( s \right){{\bf{e}}_2}\left( s \right) \\
 \,\,\,\,\,\,\,\,\,\,\,\,\,\,\,\,\,\,\,\,\,\,\,\,\,\,\,\,\,\,\,\,\,\,\,\,\,\,\,\,\,\,\,\,\,\, + \alpha \left( s \right){k_2}\left( s \right){{\bf{e}}_3}\left( s \right) \\
 \end{array}
\end{equation}
is obtained.

\noindent On the other hand, since the first normal line at the
each point of $C$ is lying in the plane generated by the second
normal line and the third normal line of ${C^*}$ at the
corresponding points under bijection $\phi $, the vector field
${{\bf{e}}_2}\left( s \right)$ is given by

$$
\begin{array}{l}
{{\bf{e}}_2}\left( s \right) = g\left( s \right){\bf{e}}_3^*\left(
{f\left( s \right)} \right) + h\left( s \right){\bf{e}}_4^*\left(
{f\left( s \right)} \right)
 \end{array}
$$
where $g$ and $h$ are some smooth functions on $L$. If we take
into consideration
$$
\begin{array}{l}
\left\langle {{\bf{e}}_1^*\left( {f\left( s \right)}
\right),\,g\left( s \right){\bf{e}}_3^*\left( {f\left( s \right)}
\right) + h\left( s \right){\bf{e}}_4^*\left( {f\left( s \right)}
\right)} \right\rangle  = 0
 \end{array}
$$
and the equation (3.4), then we have $\alpha '\left( s \right)
=0$. So we rewrite the equation (3.4) as

\begin{equation}\label{3.5}
\begin{array}{l}
f'\left( s \right){\bf{e}}_1^*\left( {f\left( s \right)} \right) =
\left( {1 + {\mu _1}\alpha {k_1}\left( s \right)}
\right){{\bf{e}}_1}\left( s \right) + \alpha {k_2}\left( s
\right){{\bf{e}}_3}\left( s \right),
 \end{array}
\end{equation}
that is,
$$
\begin{array}{l}
{\bf{e}}_1^*\left( {f\left( s \right)} \right) = \frac{{\left( {1
+ {\mu _1}\alpha {k_1}\left( s \right)} \right)}}{{f'\left( s
\right)}}{{\bf{e}}_1}\left( s \right) + \frac{{\alpha {k_2}\left(
s \right)}}{{f'\left( s \right)}}{{\bf{e}}_3}\left( s \right)
 \end{array}
$$
where
$$
\begin{array}{l}
f'\left( s \right) = \sqrt {\left| {{{\left( {1 + {\mu _1}\alpha
{k_1}\left( s \right)} \right)}^2} + \varepsilon_3 {{\left(
{\alpha {k_2}\left( s \right)} \right)}^2}} \right|}
\,,\,\,\varepsilon_3  = \left\{ \begin{array}{l}
  - 1\,\,,\,\,{{\bf{e}}_3}\,\,{\rm{is\,\,  timelike}}{\rm{,}} \\
 \,\,\,\,1\,\,\,,\,\,{{\bf{e}}_3}\,\,{\rm{is\,\,  spacelike}}{\rm{.}} \\
 \end{array} \right.
 \end{array}
$$
By taking differentiation both sides of the equations (3.5) with
respect to $s$,
\begin{equation}\label{3.6}
\begin{array}{l}
 f'\left( s \right)k_1^*\left( {f\left( s \right)} \right){\bf{e}}_2^*\left( {f\left( s \right)} \right) = {\left( {\frac{{1 + {\mu _1}\alpha {k_1}\left( s \right)}}{{f'\left( s \right)}}} \right)^\prime }{{\bf{e}}_1}\left( s \right) \\
 \,\,\,\,\,\,\,\,\,\,\,\,\,\,\,\,\,\,\,\,\,\,\,\,\,\,\,\,\,\,\,\,\,\,\,\,\,\,\,\,\,\,\,\,\,\,\,\,\,\,\,\,\,\,\,\,\,\,\,\, + \left( {\frac{{\left( {1 + {\mu _1}\alpha {k_1}\left( s \right)} \right){k_1}\left( s \right) + {\mu _2}\alpha {{\left( {{k_2}\left( s \right)} \right)}^2}}}{{f'\left( s \right)}}} \right){{\bf{e}}_2}\left( s \right) \\
 \,\,\,\,\,\,\,\,\,\,\,\,\,\,\,\,\,\,\,\,\,\,\,\,\,\,\,\,\,\,\,\,\,\,\,\,\,\,\,\,\,\,\,\,\,\,\,\,\,\,\,\,\,\,\,\,\,\,\,\,\, + {\left( {\frac{{\alpha {k_2}\left( s \right)}}{{f'\left( s \right)}}} \right)^\prime }{{\bf{e}}_3}\left( s \right) + \left( {\frac{{\alpha {k_2}\left( s \right){k_3}\left( s \right)}}{{f'\left( s \right)}}} \right){{\bf{e}}_4}\left( s \right) \\
 \end{array}
\end{equation}
is obtained for $s \in L$. Since
$$
\begin{array}{l}
\left\langle {{\bf{e}}_2^*\left( {f\left( s \right)}
\right),\,g\left( s \right){\bf{e}}_3^*\left( {f\left( s \right)}
\right) + h\left( s \right){\bf{e}}_4^*\left( {f\left( s \right)}
\right)} \right\rangle  = 0,
 \end{array}
$$
then in the equation (3.6) the coefficient of ${{\bf{e}}_2}\left(
s \right)$ vanishes, that is,
$$\begin{array}{l}
\left( {1 + {\mu _1}\alpha {k_1}\left( s \right)}
\right){k_1}\left( s \right) + {\mu _2}\alpha {\left( {{k_2}\left(
s \right)} \right)^2} = 0.
 \end{array}
$$
Thus, ${k_1}\left( s \right) =  - \alpha \left( {{\mu
_1}k_1^2\left( s \right) + {\mu _2}k_2^2\left( s \right)} \right)$
is satisfied. This completes the proof.

\noindent If we investigate the special cases separately, then we
have

in the Case 1;
$$
\begin{array}{l}
{k_1}\left( s \right) = \alpha \left( {k_1^2\left( s \right) +
k_2^2\left( s \right)} \right),
 \end{array}
$$

in the Case 2;
$$
\begin{array}{l}
{k_1}\left( s \right) = \alpha \left( {k_1^2\left( s \right) -
k_2^2\left( s \right)} \right),
 \end{array}
$$

in the Case 3;
$$
\begin{array}{l}
{k_1}\left( s \right) =  - \alpha \left( {k_1^2\left( s \right) +
k_2^2\left( s \right)} \right).
 \end{array}
$$

\begin{theorem}\label{T:3.2}
Let $C$ be a special spacelike Frenet curve in $E_1^4$ whose
curvature functions ${k_1}$ and ${k_2}$ are non-constant functions
and satisfy the equality ${k_1}\left( s \right) =  - \alpha \left(
{{\mu _1}k_1^2\left( s \right) + {\mu _2}k_2^2\left( s \right)}
\right)$, where $\alpha$ is non-zero constant, for all $s \in L$.
If the spacelike curve ${C^*}$ given by
$$
\begin{array}{l}
{{\bf{c}}^*}\left( s \right) = {\bf{c}}\left( s \right) + \alpha
{{\bf{e}}_2}\left( s \right)
 \end{array}
$$
is a special spacelike Frenet curve, then ${C^*}$ is a generalized
spacelike Mannheim mate curve of $C$.
\end{theorem}
\textbf{Proof .} The arc-length parameter of ${C^*}$ is defined by
$$
\begin{array}{l}
{s^*} = \int\limits_0^s {\left\| {\frac{{d{{\bf{c}}^*}\left( s
\right)}}{{ds}}} \right\|} ds
 \end{array}
$$
for all $s \in L$. Under the assumptation of
$$
\begin{array}{l}
{k_1}\left( s \right) =  - \alpha \left( {{\mu _1}k_1^2\left( s
\right) + {\mu _2}k_2^2\left( s \right)} \right)
 \end{array}
$$
and after calculations for all cases, separately, we obtain

in the Case 1; $\,\,\,\,f'\left( s \right) = \sqrt {\left| {1 -
\alpha {k_1}\left( s \right)} \right|}, $

in the Case 2; $\,\,\,\,f'\left( s \right) = \sqrt {\left| {1 -
\alpha {k_1}\left( s \right)} \right|}, $

in the Case 3; $\,\,\,\,f'\left( s \right) = \sqrt {\left| {1 +
\alpha {k_1}\left( s \right)} \right|}. $

\noindent Thus, we can generalize
$$
\begin{array}{l}
f'\left( s \right) = \sqrt {\left| {1 + {\mu _1}\alpha {k_1}\left(
s \right)} \right|}
 \end{array}
$$
for all $s \in L$.

\noindent By differentiating the equation ${{\bf{c}}^*}\left(
{f\left( s \right)} \right) = {\bf{c}}\left( s \right) + \alpha
{{\bf{e}}_2}\left( s \right)$ with respect to $s$, it is seen that
$$
\begin{array}{l}
f'\left( s \right){\bf{e}}_1^*\left( {f\left( s \right)} \right) =
\left( {1 + {\mu _1}\alpha {k_1}\left( s \right)}
\right){{\bf{e}}_1}\left( s \right) + \alpha {k_2}\left( s
\right){{\bf{e}}_3}\left( s \right).
 \end{array}
$$
So, it is seen that
\begin{equation}\label{3.7}
\begin{array}{l}
{\bf{e}}_1^*\left( {f\left( s \right)} \right) = \left( {\frac{{1
+ {\mu _1}\alpha {k_1}\left( s \right)}}{{\sqrt {\left| {1 + {\mu
_1}\alpha {k_1}\left( s \right)} \right|} }}{{\bf{e}}_1}\left( s
\right) + \frac{{\alpha {k_2}\left( s \right)}}{{\sqrt {\left| {1
+ {\mu _1}\alpha {k_1}\left( s \right)} \right|}
}}{{\bf{e}}_3}\left( s \right)} \right)
 \end{array}
\end{equation}
for $s \in L$.

\noindent The differentiation of the last equation with respect to
$s$ is
\begin{equation}\label{3.8}
\begin{array}{l}
 f'\left( s \right)k_1^*\left( {f\left( s \right)} \right){\bf{e}}_2^*\left( {f\left( s \right)} \right) = {\left( {\sqrt {\left| {1 + {\mu _1}\alpha {k_1}\left( s \right)} \right|} } \right)^\prime }{{\bf{e}}_1}\left( s \right) \\
 \,\,\,\,\,\,\,\,\,\,\,\,\,\,\,\,\,\,\,\,\,\,\,\,\,\,\,\,\,\,\,\,\,\,\,\,\,\,\,\,\,\,\,\, + \left( {\frac{{\left( {1 + {\mu _1}\alpha {k_1}\left( s \right)} \right){k_1}\left( s \right) + {\mu _2}\alpha k_2^2\left( s \right)}}{{\sqrt {\left| {1 + {\mu _1}\alpha {k_1}\left( s \right)} \right|} }}} \right){{\bf{e}}_2}\left( s \right) \\
 \,\,\,\,\,\,\,\,\,\,\,\,\,\,\,\,\,\,\,\,\,\,\,\,\,\,\,\,\,\,\,\,\,\,\,\,\,\,\,\,\,\,\,\,\, + {\left( {\frac{{\alpha {k_2}\left( s \right)}}{{\sqrt {\left| {1 + {\mu _1}\alpha {k_1}\left( s \right)} \right|} }}} \right)^\prime }{{\bf{e}}_3}\left( s \right) + \left( {\frac{{\alpha {k_2}\left( s \right){k_3}\left( s \right)}}{{\sqrt {\left| {1 + {\mu _1}\alpha {k_1}\left( s \right)} \right|} }}} \right){{\bf{e}}_4}\left( s \right). \\
 \end{array}
\end{equation}
According to our assumption,
$$
\begin{array}{l}
\frac{{\left( {1 + {\mu _1}\alpha {k_1}\left( s \right)}
\right){k_1}\left( s \right) + {\mu _2}\alpha k_2^2\left( s
\right)}}{{\sqrt {\left| {1 + {\mu _1}\alpha {k_1}\left( s
\right)} \right|} }} = 0
 \end{array}
$$
is hold. Thus, the coefficient of ${{\bf{e}}_2}\left( s \right)$
in the equation (3.8) is zero. It is seen from the equation (3.8),
${\bf{e}}_2^*\left( {f\left( s \right)} \right)$ is given by
linear combination of ${{\bf{e}}_1}\left( s
\right),\;\,{{\bf{e}}_3}\left( s \right)$ and ${{\bf{e}}_4}\left(
s \right)$. Also, from equation (3.7), ${\bf{e}}_1^*\left(
{f\left( s \right)} \right)$ is a linear combination of
${{\bf{e}}_1}\left( s \right)$ and ${{\bf{e}}_3}\left( s \right).$
Moreover, ${C^*}$ is a special spacelike Frenet curve that the
vector ${{\bf{e}}_2}\left( s \right)$ is given by linear
combination of ${\bf{e}}_3^*\left( {f\left( s \right)} \right)$
and ${\bf{e}}_4^*\left( {f\left( s \right)} \right)$.

\noindent Therefore, the first normal line $C$ lies in the plane
generated by the second normal line and third normal line of
${C^*}$ at the corresponding points under a bijection $\phi $
which is defined by $\phi \left( {{\bf{c}}\left( s \right)}
\right) = {{\bf{c}}^*}\left( {f\left( s \right)} \right)$. Thus,
the proof of the theorem is completed.
\begin{remark}
In 4-dimensional Minkowski space for a special spacelike Frenet
curve $C$ with curvature functions ${k_1}$ and ${k_2}$ satisfying
$$
\begin{array}{l}
{k_1}\left( s \right) =  - \alpha \left( {{\mu _1}k_1^2\left( s
\right) + {\mu _2}{\mu _3}k_2^2\left( s \right)} \right),
 \end{array}
$$
it is not clear that a smooth spacelike curve ${C^*}$ given by
(3.1) is a special Frenet curve. So, it is unknown whether the
reverse of Theorem 3.1 is true or false.
\end{remark}
\begin{theorem}
Let $C$ be a spacelike special curve in $E_1^4$ with non-zero
third curvature function ${k_3}$. If there exists a spacelike
special Frenet curve ${C^*}$ in $E_1^4$ such that the first normal
line of $C$ is linearly dependent with the third normal line of
${C^*}$ at the corresponding points $\bf{c}\left( s \right)$ and
${\bf{c}^*}\left( s \right)$, respectively, under a bijection
$\phi :C \to {C^*}$, then the curvatures ${k_1}$ and ${k_2}$ of
$C$ are constant functions.
\end{theorem}
\textbf{Proof.} Let $C$ be a spacelike Frenet curve in $E_1^4$
with the Frenet frame field $\left\{
{{{\bf{e}}_1},\,{{\bf{e}}_2},\,{{\bf{e}}_3},\,{{\bf{e}}_4}}
\right\}$ and curvature functions ${k_1},\,{k_2}$ and ${k_3}$.
Also, we assume that ${C^*}$ be a spacelike special Frenet curve
in $E_1^4$ with the Frenet frame field $\left\{
{{\bf{e}}_1^*,\,{\bf{e}}_2^*,\,{\bf{e}}_3^*,\,{\bf{e}}_4^*}
\right\}$ and curvature functions $k_1^*,\,k_2^*\,$ and $k_3^*$.

\noindent Let the first normal line of $C$ be linearly dependent
with the third normal line of ${C^*}$ at the corresponding points
$C$ and ${C^*}$, respectively. Then the parametrization of ${C^*}$
is
\begin{equation}\label{3.9}
\begin{array}{l}
{{\bf{c}}^*}\left( {f\left( s \right)} \right) = {\bf{c}}\left( s
\right) + \alpha \left( s \right){{\bf{e}}_2}\left( s \right)
 \end{array}
\end{equation}
for all $s \in L$. If ${s^*}$ is the arc-length parameter of
${C^*}$, then
\begin{equation}\label{3.10}
\begin{array}{l}
{s^*} = \int\limits_0^s {\sqrt {\left| {{{\left( {1 + {\mu
_1}\alpha {k_1}} \right)}^2} + \varepsilon_2 \left( {\alpha
'\left( s \right)} \right) + \varepsilon_3 {{\left( {\alpha \left(
s \right){k_2}\left( s \right)} \right)}^2}} \right|} } ds
 \end{array}
\end{equation}
where
$$
\begin{array}{l}
\varepsilon_i  = \left\{ \begin{array}{l}
  - 1\,\,,\,\,{{\bf{e}}_i}\,\,{\rm{is}}\,\,{\rm{timelike}} \\
 \,\,\,1\,\,\,,\,\,{{\bf{e}}_i}\,\,{\rm{is}}\,\,{\rm{spacelike}} \\
 \end{array} \right.,\,\,\,{\rm{for}}\,\,\,\,\, {\rm{ }}2 \le i \le 4
 \end{array}
$$
and
$$
\begin{array}{l}
f:\,L \to {L^*} \\
 \,\,\,\,\,\,\,s\,\, \to \,f\left( s \right) = {s^*}. \\
 \end{array}
$$
Moreover, $\phi :C \to {C^*}$ is a bijection given by $\phi \left(
{{\bf{c}}\left( s \right)} \right) = {{\bf{c}}^*}\left( {f\left( s
\right)} \right)$.

\noindent By differentiating the equation (3.9) with respect to
$s$ and using Frenet formulas, we have
\begin{equation}\label{3.11}
\begin{array}{l}
 f'\left( s \right){\bf{e}}_1^*\left( {f\left( s \right)} \right) = \left( {1 + {\mu _1}\alpha \left( s \right){k_1}\left( s \right)} \right){{\bf{e}}_1}\left( s \right) + \alpha '\left( s \right){{\bf{e}}_2}\left( s \right) \\
 \,\,\,\,\,\,\,\,\,\,\,\,\,\,\,\,\,\,\,\,\,\,\,\,\,\,\,\,\,\,\,\,\,\,\,\,\,\,\,\,\,\,\,\,\,\, + \alpha \left( s \right){k_2}\left( s \right){{\bf{e}}_3}\left( s \right). \\
 \end{array}
\end{equation}
Since ${\bf{e}}_4^*\left( {f\left( s \right)} \right) =  \mp
{{\bf{e}}_2}\left( s \right)$, then
$$
\begin{array}{l}
 \left\langle {f'\left( s \right){\bf{e}}_1^*\left( {f\left( s \right)} \right),\,{\bf{e}}_4^*\left( {f\left( s \right)} \right)} \right\rangle  = \left\langle {\left( {1 + {\mu _1}\alpha \left( s \right){k_1}\left( s \right)} \right){{\bf{e}}_1}\left( s \right) + \alpha '\left( s \right){{\bf{e}}_2}\left( s \right)} \right. \\
 \,\,\,\,\,\,\,\,\,\,\,\,\,\,\,\,\,\,\,\,\,\,\,\,\,\,\,\,\,\,\,\,\,\,\,\,\,\,\,\,\,\,\,\,\,\,\,\,\,\,\,\,\,\,\,\,\,\,\,\,\,\,\,\,\,\,\,\,\,\,\,\,\,\,\,\,\,\,\,\,\,\,\,\left. { + \alpha \left( s \right){k_2}\left( s \right){{\bf{e}}_3}\left( s \right),\, \mp {{\bf{e}}_2}\left( s \right)} \right\rangle , \\
 \end{array}
$$
that is,
$$
\begin{array}{l}
0 =  \mp \alpha '\left( s \right).
 \end{array}
$$
It is easily seen that $\alpha $ is a constant number from the
last equation. Thus, hereafter we can denote $\alpha \left( s
\right) = \alpha $, for all $s \in L.$

\noindent From the equation (3.10), we get
$$
\begin{array}{l}
f'\left( s \right) = \sqrt {\left| {{{\left( {1 + {\mu _1}\alpha
{k_1}\left( s \right)} \right)}^2} + \varepsilon_3 {{\left(
{\alpha {k_2}\left( s \right)} \right)}^2}} \right|}  > 0
 \end{array}
$$
where
$$
\begin{array}{l}
\varepsilon_3 = \left\{ \begin{array}{l}
 - 1\,\,,\,\,{{\bf{e}}_i}\,\,{\rm{is}}\,\,{\rm{timelike}} \\
 \,\,\,\,1\,\,,\,\,{{\bf{e}}_i}\,\,{\rm{is}}\,\,{\rm{spacelike}} \\
 \end{array} \right. \,\,\,\,,\,{\rm{for   }}\,\,\,\,\, 2 \le i \le 4.
 \end{array}
$$
Then, we rewrite the equation (3.11) as follows;
$$
\begin{array}{l}
{\bf{e}}_1^*\left( {f\left( s \right)} \right) = \left( {\frac{{1
+ {\mu _1}\alpha {k_1}\left( s \right)}}{{f'\left( s \right)}}}
\right){{\bf{e}}_1}\left( s \right) + \left( {\frac{{\alpha
{k_2}\left( s \right)}}{{f'\left( s \right)}}}
\right){{\bf{e}}_3}\left( s \right).\
 \end{array}
$$
The differentiation of the last equation with respect to $s$ is
\begin{equation}\label{3.12}
\begin{array}{l}
 f'\left( s \right)k_1^*\left( {f\left( s \right)} \right){\bf{e}}_2^*\left( {f\left( s \right)} \right) = {\left( {\frac{{1 + {\mu _1}\alpha {k_1}\left( s \right)}}{{f'\left( s \right)}}} \right)^\prime }{{\bf{e}}_1}\left( s \right) \\
 \,\,\,\,\,\,\,\,\,\,\,\,\,\,\,\,\,\,\,\,\,\,\,\,\,\,\,\,\,\,\,\,\,\,\,\,\,\,\,\,\,\,\,\,\,\,\,\,\,\,\,\,\,\,\,\,\,\,\,\,\,\,\,\,\,\,\,\,\, + \left( {\frac{{{k_1}\left( s \right) + {\mu _1}\alpha k_1^2\left( s \right) + {\mu _2}\alpha k_2^2\left( s \right)}}{{f'\left( s \right)}}} \right){{\bf{e}}_2}\left( s \right) \\
 \,\,\,\,\,\,\,\,\,\,\,\,\,\,\,\,\,\,\,\,\,\,\,\,\,\,\,\,\,\,\,\,\,\,\,\,\,\,\,\,\,\,\,\,\,\,\,\,\,\,\,\,\,\,\,\, \,\,\,\,\,\,\,\,\,\,\,\,\,+ {\left( {\frac{{\alpha {k_2}\left( s \right)}}{{f'\left( s \right)}}} \right)^\prime }{{\bf{e}}_3}\left( s \right) + \left( {\frac{{\alpha {k_2}\left( s \right){k_3}\left( s \right)}}{{f'\left( s \right)}}} \right){{\bf{e}}_4}\left( s \right). \\
 \end{array}
\end{equation}
Since $\left\langle {f'\left( s \right)k_1^*\left( {f\left( s
\right)} \right){\bf{e}}_2^*\left( {f\left( s \right)}
\right),\,{\bf{e}}_4^*\left( {f\left( s \right)} \right)}
\right\rangle  = 0$ and ${\bf{e}}_4^*\left( {f\left( s \right)}
\right) =  \mp {{\bf{e}}_2}\left( s \right)$ for all $s \in L$, we
obtain
$$
\begin{array}{l}
{k_1}\left( s \right) + {\mu _1}\alpha k_1^2\left( s \right) +
{\mu _2}\alpha k_2^2\left( s \right) = 0
 \end{array}
$$
is satisfied. Then,
\begin{equation}\label{3.13}
\begin{array}{l}
\alpha  =  - \frac{{{k_1}\left( s \right)}}{{{\mu _1} k_1^2\left(
s \right) + {\mu _2} k_2^2\left( s \right)}}
 \end{array}
\end{equation}
is a non-zero constant number. Thus, from the equation (3.12), it
is seen that
$$
\begin{array}{l}
 {\bf{e}}_2^*\left( {f\left( s \right)} \right) = \frac{1}{{f'\left( s \right)K\left( s \right)}}{\left( {\frac{{1 + {\mu _1}\alpha {k_1}\left( s \right)}}{{f'\left( s \right)}}} \right)^\prime }{{\bf{e}}_1}\left( s \right) + \frac{1}{{f'\left( s \right)K\left( s \right)}}\left( {\frac{{\alpha {k_2}\left( s \right)}}{{f'\left( s \right)}}} \right){{\bf{e}}_3}\left( s \right) \\
 \,\,\,\,\,\,\,\,\,\,\,\,\,\,\,\,\,\,\,\,\,\,\,\,\,\,\,\,\,\, + \frac{1}{{f'\left( s \right)K\left( s \right)}}\left( {\frac{{\alpha {k_2}\left( s \right){k_3}\left( s \right)}}{{f'\left( s \right)}}} \right){{\bf{e}}_4}\left( s \right) \\
\end{array}
$$
where $K\left( s \right) = k_1^*\left( {f\left( s \right)}
\right)$ for all $s \in L$. By differentiating the last equation
with respect to $s$, we obtain
$$
\begin{array}{l}
 f'\left( s \right)\left[ {{\mu _1}k_1^*\left( {f\left( s \right)} \right){\bf{e}}_1^*\left( {f\left( s \right)} \right) + k_2^*\left( {f\left( s \right)} \right){{\bf{e}}_3^*}\left( {f\left( s \right)} \right)} \right] = {\left( {\frac{1}{{f'\left( s \right)K\left( s \right)}}{{\left( {\frac{{1 + {\mu _1}\alpha {k_1}\left( s \right)}}{{f'\left( s \right)}}} \right)}^\prime }} \right)^\prime }{{\bf{e}}_1}\left( s \right) \\
 \,\,\,\,\,\,\,\,\,\,\,\,\,\,\,\,\,\,\,\,\,\,\,\,\,\,\,\,\,\,\,\,\,\,\,\,\,\,\,\,\,\,\,\,\,\,\,\,\, + \left( {\frac{k_1{\left(s \right)}}{{f'\left( s \right)K\left( s \right)}}{{\left( {\frac{{1 + {\mu _1}\alpha {k_1}\left( s \right)}}{{f'\left( s \right)}}} \right)}^\prime } + \frac{{{\mu _2}{k_2}\left( s \right)}}{{f'\left( s \right)K\left( s \right)}}{{\left( {\frac{{\alpha {k_2}\left( s \right)}}{{f'\left( s \right)}}} \right)}^\prime }} \right){{\bf{e}}_2}\left( s \right) \\
 \,\,\,\,\,\,\,\,\,\,\,\,\,\,\,\,\,\,\,\,\,\,\,\,\,\,\,\,\,\,\,\,\,\,\,\,\,\,\,\,\,\,\,\,\,\,\,\,\, + \left( {{{\left( {\frac{1}{{f'\left( s \right)K\left( s \right)}}{{\left( {\frac{{\alpha {k_2}\left( s \right)}}{{f'\left( s \right)}}} \right)}^\prime }} \right)}^\prime } + \frac{{{\mu _3}{k_3}\left( s \right)}}{{f'\left( s \right)K\left( s \right)}}\left( {\frac{{\alpha {k_2}\left( s \right){k_3}\left( s \right)}}{{f'\left( s \right)}}} \right)} \right){{\bf{e}}_3}\left( s \right) \\
 \,\,\,\,\,\,\,\,\,\,\,\,\,\,\,\,\,\,\,\,\,\,\,\,\,\,\,\,\,\,\,\,\,\,\,\,\,\,\,\,\,\,\,\,\,\,\,\,\, + \left( {{{\left( {\frac{1}{{f'\left( s \right)K\left( s \right)}}{{\left( {\frac{{\alpha {k_2}\left( s \right){k_3}\left( s \right)}}{{f'\left( s \right)}}} \right)}^\prime }} \right)}^\prime } + \frac{{{k_3}\left( s \right)}}{{f'\left( s \right)K\left( s \right)}}{{\left( {\frac{{\alpha {k_2}\left( s \right)}}{{f'\left( s \right)}}} \right)}^\prime }} \right){{\bf{e}}_4}\left( s \right) \\
\end{array}
$$
for all $s \in L$. If we take into consideration
$$
\begin{array}{l}
\left\langle {f'\left( s \right)\left( {{\mu _1}k_1^*\left(
{f\left( s \right)} \right){\bf{e}}_1^*\left( {f\left( s \right)}
\right) + k_2^*\left( {f\left( s \right)}
\right){{\bf{e}}_3^*}\left( {f\left( s \right)} \right)}
\right)\,,\,\,{\bf{e}}_4^*\left( {f\left( s \right)} \right)}
\right\rangle  = 0\
\end{array}
$$
and
$$
\begin{array}{l}
{\bf{e}}_4^*\left( {f\left( s \right)} \right) =  \mp
{{\bf{e}}_2}\left( s \right),\
\end{array}
$$
then
$$
\begin{array}{l}
 {\mu _1}\alpha {k_1}\left( s \right){k_1}^\prime \left( s \right) {f'\left( s \right)}- {k_1}\left( s \right)\left( {1 + {\mu _1}\alpha {k_1}\left( s \right)} \right)f''\left( s \right) \\
 \,\,\,\,\,\,\,\,\,\,\,\,\,\,\,\,\,\,\,\,\,\,\,\,\,\,\,\,\,\,\,\,\,\,\,\,\,\,\,\,\,\,\,\,\,\,\,\, + {\mu _2}\alpha {k_2}\left( s \right){k_2}^\prime \left( s \right){f'\left( s \right)} - {\mu _2}\alpha k_2^2\left( s \right)f''\left( s \right) = 0. \\
\end{array}
$$
If we arrange the last equation, then we find
\begin{equation}\label{3.14}
\begin{array}{l}
 \alpha \left( {{\mu _1}{k_1}\left( s \right){{k'}_1}\left( s \right) + {\mu _2}{k_2}\left( s \right){{k'}_2}\left( s \right)} \right)f'\left( s \right) \\
 \,\,\,\,\,\,\,\,\,\,\,\,\,\,\,\,\,\,\,\,\, - \left( {{k_1} + \alpha \left( {{\mu _1}k_1^2\left( s \right) + {\mu _2}k_2^2\left( s \right)} \right)} \right)f''\left( s \right) = 0. \\
 \end{array}
\end{equation}
Moreover, the differentiation of the equation (3.13) with respect
to $s$ is
$$
\begin{array}{l}
{k'_1}\left( s \right) + 2\alpha \left( {{\mu _1}{k_1}\left( s
\right){{k'}_1}\left( s \right) + {\mu _2}{k_2}\left( s
\right){{k'}_2}\left( s \right)} \right) = 0.
\end{array}
$$
From the above equation, we see
\begin{equation}\label{3.15}
\begin{array}{l}
 - \frac{{{{k'}_1}\left( s \right)}}{2} = \alpha \left( {{\mu
_1}{k_1}\left( s \right){{k'}_1}\left( s \right) + {\mu
_2}{k_2}\left( s \right){{k'}_2}\left( s \right)} \right).
 \end{array}
\end{equation}
If we substitute the equations (3.13) and (3.15) into the equation
(3.14), we obtain
$$
\begin{array}{l}
 - \frac{{{{k'}_1}\left( s \right)}}{2} = 0.
\end{array}
$$
Finally, we find that the first curvature function is constant
(that is, positive constant).

\noindent Thus, from the equation (3.15) it is seen that the
second curvature function ${k_2}$ is positive constant, too. This
completes the proof.

\noindent In \cite{Eis}, a formula of parametric equation of
Mannheim curve is given in ${E^3}$. Moreover, the parametric
equation of generalized Mannheim curve in ${E^4}$ is obtained in
\cite{Mat2}. The following theorem gives a parametric
representation of a generalized spacelike Mannheim curve with
timelike second binormal vector in $E_1^4$.

\begin{theorem}
Let   be a spacelike special curve defined by
$$
\begin{array}{l}
{\bf{c}}\left( u \right) = \left[ {\begin{array}{*{20}{c}}
   {\alpha \int {f\left( u \right)\sinh udu} }  \\
   {\alpha \int {f\left( u \right)\cosh udu} }  \\
   {\alpha \int {f\left( u \right)g\left( u \right)du} }  \\
   {\alpha \int {f\left( u \right)h\left( u \right)du} }  \\
\end{array}} \right]\
\end{array}
$$
for $u \in I \subset\mathbb{R}$. Here $\alpha$ is a non-zero
constant number, $g:I \to \mathbb{R}$ and $h:I \to \mathbb{R}$ are
any smooth functions and the positive valued smooth function $f:I
\to \mathbb{R}$ is given by

$$
\begin{array}{l}
 f\left( u \right) = \left( {1 + g^2 \left( u \right) + h^2 \left( u \right)} \right)^{{{ - 3} \mathord{\left/
 {\vphantom {{ - 3} 2}} \right.
 \kern-\nulldelimiterspace} 2}}  \\
 \,\,\,\,\,\,\,\,\,\,\,\,\,\,\,\,\,\,\left| { - 1 - g^2 \left( u \right) - h^2 \left( u \right) + \dot g^2 \left( u \right) + \dot h^2 \left( u \right) + \left( {\dot g\left( u \right)h\left( u \right) - g\left( u \right)\dot h\left( u \right)} \right)^2 } \right|^{{{ - 5} \mathord{\left/
 {\vphantom {{ - 5} 2}} \right.
 \kern-\nulldelimiterspace} 2}}  \\
 \,\,\,\,\,\,\,\,\,\,\,\,\,\,\,\,\,\;\left| {\left( { - 1 - g^2 \left( u \right) - h^2 \left( u \right) + \dot g^2 \left( u \right) + \dot h^2 \left( u \right) + \left( {\dot g\left( u \right)h\left( u \right) - g\left( u \right)\dot h\left( u \right)} \right)^2 } \right)^3 } \right. \\
 \quad \,\,\,\,\,\,\,\,\,\,\,\,\,\,\, - \left( {1 + g^2 \left( u \right) + h^2 \left( u \right)} \right)^3 \left[ {\left( {g\left( u \right) - \ddot g\left( u \right)} \right)^2  + \left( {h\left( u \right) - \ddot h\left( u \right)} \right)^2 } \right. \\
 \,\,\,\,\,\,\,\,\,\,\,\,\,\,\,\,\,\,\,\left. {\, - \left( {\left( {g\left( u \right)\dot h\left( u \right) - \dot g\left( u \right)h\left( u \right)} \right) + \left( {\dot g\left( u \right)\ddot h\left( u \right) - \ddot g\left( u \right)\dot h\left( u \right)} \right)} \right)^2  + \left( {g\left( u \right)\ddot h\left( u \right) - \ddot g\left( u \right)h} \right)\left( u \right)^2 } \right| \\
 \quad \quad \quad \quad \quad \quad \;\;\quad \quad \quad \;\,\,\,\,\,\,\,\,\,\,\,\,\,\,\,\,\,\,\,\,\,\,\,\,\,\,\,\,\,\,\,\,\,\,\,\,\,\,, \\
 \end{array}
$$ for $u \in I$. Then the curvature functions
${k_1}$ and ${k_2}$ of $C$ satisfy

$$
\begin{array}{l}
{k_1}\left( u \right) = \alpha \left( {k_1^2\left( u \right) +
k_2^2\left( u \right)} \right)
\end{array}
$$
at the each point ${\bf{c}}\left( u \right)$ of $C$.
\end{theorem}
\textbf{Proof.} Let $C$ be a spacelike special curve defined by
$$
\begin{array}{l}
{\bf{c}}\left( u \right) = \left[ {\begin{array}{*{20}{c}}
   {\alpha \int {f\left( u \right)\sinh udu} }  \\
   {\alpha \int {f\left( u \right)\cosh udu} }  \\
   {\alpha \int {f\left( u \right)g\left( u \right)du} }  \\
   {\alpha \int {f\left( u \right)h\left( u \right)du} }  \\
\end{array}} \right]\quad ,\quad u \in I \subset\mathbb{R}
\end{array}
$$
where $\alpha$ is a non-zero constant number, $g$ and $ h$ are any
smooth functions. $f$ is a positive valued smooth function. Thus,
we obtain
\begin{equation}\label{3.16}
\begin{array}{l}
{\bf{\dot c}}\left( u \right) = \left[ {\begin{array}{*{20}{c}}
   {\alpha f\left( u \right)\sinh u}  \\
   {\alpha f\left( u \right)\cosh u}  \\
   {\alpha f\left( u \right)g\left( u \right)}  \\
   {\alpha f\left( u \right)h\left( u \right)}  \\
\end{array}} \right]\quad ,\quad u \in I \subset\mathbb{R}
 \end{array}
\end{equation}
where the subscript dot (.) denotes the differentiation with
respect to $u$.

\noindent The arc-length parameter $s$ of $C$ is given by
$$
\begin{array}{l}
s = \psi \left( u \right) = \int\limits_{{u_0}}^u {\left\|
{{\bf{\dot c}}\left( u \right)} \right\|} du\
\end{array}
$$
where $\left\| {{\bf{\dot c}}\left( u \right)} \right\| = \alpha
f\left( u \right)\sqrt {1 + {g^2}\left( u \right) + {h^2}\left( u
\right)}.$

\noindent If $\varphi$ denotes the inverse function of $\psi :I
\to L \subset\mathbb{R}$, then $u = \varphi \left( s \right)$ and
$$
\begin{array}{l}
\varphi '\left( s \right) = {\left\| {{{\left.
{\frac{{d{\bf{c}}\left( u \right)}}{{du}}} \right|}_{u = \varphi
\left( s \right)}}} \right\|^{ - 1}}\quad ,\quad s \in I
\end{array}
$$
where the prime $\left( ' \right)$ denotes the differentiation
with respect to $s$.

\noindent The unit tangent vector ${{\bf{e}}_1}\left( s \right)$
of the curve $C$ at the each point ${\bf{c}}\left( {\varphi \left(
s \right)} \right)$ is given by
\begin{equation}\label{3.17}
\begin{array}{l}
{{\bf{e}}_1}\left( s \right) = {\left( {1 + {g^2}\left( {\varphi
\left( s \right)} \right) + {h^2}\left( {\varphi \left( s \right)}
\right)} \right)^{ - {1 \mathord{\left/
 {\vphantom {1 2}} \right.
 \kern-\nulldelimiterspace} 2}}}\left[ {\begin{array}{*{20}{c}}
   {\sinh \left( {\varphi \left( s \right)} \right)}  \\
   {\cosh \left( {\varphi \left( s \right)} \right)}  \\
   {g\left( {\varphi \left( s \right)} \right)}  \\
   {h\left( {\varphi \left( s \right)} \right)}  \\
\end{array}} \right]
\end{array}
\end{equation}
for all $s \in L$. Some simplifying assumptions are made for the
sake of brevity as follows;
$$
\begin{array}{l}
 \sinh : = \sinh \left( {\varphi \left( s \right)} \right)\quad ,\quad \quad  \cosh : = \cosh \left( {\varphi \left( s \right)} \right) \\
 f: = f\left( {\varphi \left( s \right)} \right)\quad \quad \quad \,\, , \quad \quad g: = g\left( {\varphi \left( s \right)} \right)\quad \quad \quad ,\quad h: = h\left( {\varphi \left( s \right)} \right), \\
 \dot g: = \dot g\left( {\varphi \left( s \right)} \right) = {\left. {\frac{{dg\left( u \right)}}{{du}}} \right|_{u = \varphi \left( s \right)}}\quad ,\quad \dot h: = \dot h\left( {\varphi \left( s \right)} \right) = {\left. {\frac{{dh\left( u \right)}}{{du}}} \right|_{u = \varphi \left( s \right)}}, \\
 \ddot g: = \ddot g\left( {\varphi \left( s \right)} \right) = {\left. {\frac{{{d^2}g\left( u \right)}}{{d{u^2}}}} \right|_{u = \varphi \left( s \right)}}\quad ,\quad \ddot h: = \ddot h\left( {\varphi \left( s \right)} \right) = {\left. {\frac{{{d^2}h\left( u \right)}}{{d{u^2}}}} \right|_{u = \varphi \left( s \right)}}, \\
 \varphi ': = \varphi '\left( s \right) = {\left. {\frac{{d\varphi }}{{ds}}} \right|_s}, \\
 A: = 1 + {g^2} + {h^2}\quad \,\,\,,\quad \quad B: = g\dot g + h\dot h\quad ,\quad C: = {{\dot g}^2} + {{\dot h}^2}, \\
 D: = g\ddot g + h\ddot h\quad \quad \quad ,\quad \quad E: = \dot g\ddot g + \dot h\ddot h\quad ,\quad F: = {{\ddot g}^2} + {{\ddot h}^2}.\quad  \\
\end{array}
$$
Then, we have
$$
\begin{array}{l}
\dot A = 2B\quad ,\quad \dot B = C + D\quad ,\quad \dot C =
2E\quad ,\quad \varphi ' = {\alpha ^{ - 1}}{f^{ - 1}}{A^{{{ - 1}
\mathord{\left/
 {\vphantom {{ - 1} 2}} \right.
 \kern-\nulldelimiterspace} 2}}}.
\end{array}
$$
So, we rewrite the equation (3.17) as
\begin{equation}\label{3.18}
\begin{array}{l}
{{\bf{e}}_1}: = {{\bf{e}}_1}\left( s \right) = {A^{ - {1
\mathord{\left/
 {\vphantom {1 2}} \right.
 \kern-\nulldelimiterspace} 2}}}\left[ {\begin{array}{*{20}{c}}
   {\sinh }  \\
   {\cosh }  \\
   g  \\
   h  \\
\end{array}} \right].
 \end{array}
\end{equation}
By differentiating the last equation with respect to $s$, we find
$$
\begin{array}{l}
{{\bf{e'}}_1} = \varphi '\left[ {\begin{array}{*{20}{c}}
   { - \frac{1}{2}{A^{ - {3 \mathord{\left/
 {\vphantom {3 2}} \right.
 \kern-\nulldelimiterspace} 2}}}\dot A\sinh  + {A^{ - {1 \mathord{\left/
 {\vphantom {1 2}} \right.
 \kern-\nulldelimiterspace} 2}}}\cosh }  \\
   { - \frac{1}{2}{A^{ - {3 \mathord{\left/
 {\vphantom {3 2}} \right.
 \kern-\nulldelimiterspace} 2}}}\dot A\cosh  + {A^{ - {1 \mathord{\left/
 {\vphantom {1 2}} \right.
 \kern-\nulldelimiterspace} 2}}}\sinh }  \\
   { - \frac{1}{2}{A^{ - {3 \mathord{\left/
 {\vphantom {3 2}} \right.
 \kern-\nulldelimiterspace} 2}}}\dot Ag + {A^{ - {1 \mathord{\left/
 {\vphantom {1 2}} \right.
 \kern-\nulldelimiterspace} 2}}}\dot g}  \\
   { - \frac{1}{2}{A^{ - {3 \mathord{\left/
 {\vphantom {3 2}} \right.
 \kern-\nulldelimiterspace} 2}}}\dot Ah + {A^{ - {1 \mathord{\left/
 {\vphantom {1 2}} \right.
 \kern-\nulldelimiterspace} 2}}}\dot h}  \\
\end{array}} \right],
\end{array}
$$
that is,

\begin{equation}\label{3.18}
\begin{array}{l}
{{\bf{e'}}_1} =  - \varphi '{A^{ - {1 \mathord{\left/
 {\vphantom {1 2}} \right.
 \kern-\nulldelimiterspace} 2}}}\left[ {\begin{array}{*{20}{c}}
   {{A^{ - 1}}B\sinh  - \cosh }  \\
   {{A^{ - 1}}B\cosh  - \sinh }  \\
   {{A^{ - 1}}Bg - \dot g}  \\
   {{A^{ - 1}}Bh - \dot h}  \\
\end{array}} \right].
 \end{array}
\end{equation}
From the last equation, we obtain
\begin{equation}\label{3.20}
\begin{array}{l}
{k_1}: = {k_1}\left( s \right) = \left\| {{{{\bf{e'}}}_1}\left( s
\right)} \right\| = \varphi '{A^{ - 1}}{\left| { - A + AC - {B^2}}
\right|^{{1 \mathord{\left/
 {\vphantom {1 2}} \right.
 \kern-\nulldelimiterspace} 2}}}.
 \end{array}
\end{equation}
By the fact that ${{\bf{e}}_2}\left( s \right) = {\left(
{{k_1}\left( s \right)} \right)^{ - 1}}{{\bf{e'}}_1}\left( s
\right)$, we have
$$
\begin{array}{l}
{{\bf{e}}_2}: = {{\bf{e}}_2}\left( s \right) =  - {A^{{1
\mathord{\left/
 {\vphantom {1 2}} \right.
 \kern-\nulldelimiterspace} 2}}}{\left| { - A + AC - {B^2}} \right|^{{{ - 1} \mathord{\left/
 {\vphantom {{ - 1} 2}} \right.
 \kern-\nulldelimiterspace} 2}}}\left[ {\begin{array}{*{20}{c}}
   {{A^{ - 1}}B\sinh  - \cosh }  \\
   {{A^{ - 1}}B\cosh  - \sinh }  \\
   {{A^{ - 1}}Bg - \dot g}  \\
   {{A^{ - 1}}Bh - \dot h}  \\
\end{array}} \right].
\end{array}
$$
In order to get second curvature function ${k_2}$, we need to
calculate ${k_2}\left( s \right) = \left\| {{{{\bf{e'}}}_2}\left(
s \right) - {\mu _1}{k_1}\left( s \right){{\bf{e}}_1}\left( s
\right)} \right\|$. It is seen from the above equation
$\left\langle {{{\bf{e}}_2}\left( s \right),{{\bf{e}}_2}\left( s
\right)} \right\rangle  = 1$, that is, ${{\bf{e}}_2}$ is
spacelike. Thus, ${\mu _1}$ is equal to $- 1$ and ${k_2}\left( s
\right) = \left\| {{{{\bf{e'}}}_2}\left( s \right) + {k_1}\left( s
\right){{\bf{e}}_1}\left( s \right)} \right\|$. After a long
process of calculation, we have
\begin{equation}\label{3.21}
\begin{array}{l}
{{\bf{e'}}_2} + {k_1}{{\bf{e}}_1} = \varphi '{A^{{{ - 3}
\mathord{\left/
 {\vphantom {{ - 3} 2}} \right.
 \kern-\nulldelimiterspace} 2}}}{\left| { - A + AC - {B^2}} \right|^{{{ - 3} \mathord{\left/
 {\vphantom {{ - 3} 2}} \right.
 \kern-\nulldelimiterspace} 2}}}\left[ {\begin{array}{*{20}{c}}
   {\left( {P + Q} \right)\sinh  - R\cosh }  \\
   {\left( {P + Q} \right)\cosh  - R\sinh }  \\
   {Pg - R\dot g + Q\ddot g}  \\
   {Ph - R\dot h + Q\ddot h}  \\
\end{array}} \right]
 \end{array}
\end{equation}
where
\begin{equation}\label{3.22}
\begin{array}{l}
 P = {\left( { - A + AC - {B^2}} \right)^2} + \left( { - A + AC - {B^2}} \right)\left( {{B^2} - AC - AD} \right) \\
 \,\,\,\,\,\,\,\,\,\,\,\,\, + AB\left( { - B + AE - BD} \right), \\
 Q = {A^2}\left( { - A + AC - {B^2}} \right), \\
 R = {A^2}\left( { - B + AE - BD} \right). \\
 \end{array}
\end{equation}
If we simplify $P$ then we have

$$
\begin{array}{l}
P = {A^2}\left( {1 - C + BE + D - CD} \right).
\end{array}
$$
Thus, we rewrite the equations (3.22) and (3.23) as

\begin{equation}\label{3.23}
\begin{array}{l}
{{\bf{e'}}_2} + {k_1}{{\bf{e}}_1} = \varphi '{A^{{1
\mathord{\left/
 {\vphantom {1 2}} \right.
 \kern-\nulldelimiterspace} 2}}}{\left| { - A + AC - {B^2}} \right|^{{{ - 3} \mathord{\left/
 {\vphantom {{ - 3} 2}} \right.
 \kern-\nulldelimiterspace} 2}}}\left[ {\begin{array}{*{20}{c}}
   {\left( {\tilde P + \tilde Q} \right)\sinh  - \tilde R\cosh }  \\
   {\left( {\tilde P + \tilde Q} \right)\cosh  - \tilde R\sinh }  \\
   {\tilde Pg - \tilde R\dot g + \tilde Q\ddot g}  \\
   {\tilde Ph - \tilde R\dot h + \tilde Q\ddot h}  \\
\end{array}} \right]
 \end{array}
\end{equation}
where
\begin{equation}\label{3.23}
\begin{array}{l}
 \tilde P = 1 - C + BE + D - CD, \\
 \tilde Q =  - A + AC - {B^2}, \\
 \tilde R =  - B + AE - BD. \\
 \end{array}
\end{equation}
Consequently, from the equations (3.24) and (3.25), we find

$$
\begin{array}{l}
 {\left\| {{{{\bf{e'}}}_2} + {k_1}{{\bf{e}}_1}} \right\|^2} = {\left( {\varphi '} \right)^2}A{\left| { - A + AC - {B^2}} \right|^{ - 3}}\,\left| {{{\left( {\tilde P + \tilde Q} \right)}^2} - {{\tilde R}^2}} \right. \\
 \quad \quad \quad \quad \quad \, + {{\tilde P}^2}\left( {{g^2} + {h^2}} \right) + {{\tilde R}^2}\left( {{{\dot g}^2} + {{\dot h}^2}} \right) + {{\tilde Q}^2}\left( {{{\ddot g}^2} + {{\ddot h}^2}} \right) \\
 \quad \quad \quad \quad \quad \,\left. { - 2\tilde P\tilde R\left( {g\dot g + h\dot h} \right) - 2\tilde R\tilde Q\left( {\dot g\ddot g + \dot h\ddot h} \right) + 2\tilde P\tilde Q\left( {g\ddot g + h\ddot h} \right)} \right|. \\
 \end{array}
$$
If we substitute the abbreviations into the last equation, we get

$$
\begin{array}{l}
 {\left\| {{{{\bf{e'}}}_2} + {k_1}{{\bf{e}}_1}} \right\|^2} = {\left( {\varphi '} \right)^2}A{\left| { - A + AC - {B^2}} \right|^{ - 3}} \\
 \,\,\,\,\,\,\,\,\,\,\,\,\,\,\,\,\,\,\,\,\,\,\,\,\,\,\,\left| {{{\tilde P}^2}A + 2\tilde P\tilde Q + {{\tilde Q}^2} - {{\tilde R}^2} + {{\tilde R}^2}C} \right. \\
 \,\,\,\,\,\,\,\,\,\,\,\,\,\,\,\,\,\,\,\,\,\,\,\,\,\,\,\left. { + {{\tilde Q}^2}F - 2\tilde P\tilde RB - 2\tilde R\tilde QE + 2\tilde P\tilde QD} \right|. \\
 \end{array}
$$
After substituting the equation (3.24) into the last equation and
simplifying it, we have
$$
\begin{array}{l}
 k_2^2 = {\left\| {{{{\bf{e'}}}_2} + {k_1}{{\bf{e}}_1}} \right\|^2} \\
 \quad \, = {\left( {\varphi '} \right)^2}A{\left| { - A + AC - {B^2}} \right|^{ - 2}}\,\,\left| {\left( { - A + AC - {B^2}} \right)\left( {1 + F} \right)} \right. \\
 \,\,\,\,\,\,\,\,\,\,\left. { + \left( {1 - C} \right){{\left( {1 + D} \right)}^2} + 2BE\left( {1 + D} \right) - A{E^2}} \right|\,. \\
 \end{array}
$$
Moreover, from the equation (3.20) it is seen that
$$
\begin{array}{l}
k_1^2 = {\left( {\varphi '} \right)^2}{A^{ - 2}}\left| { - A + AC
- {B^2}} \right|.
 \end{array}
$$
The last two equation gives us
$$
\begin{array}{l}
 k_1^2 + k_2^2 = {\left( {\varphi '} \right)^2}{A^{ - 2}}{\left| { - A + AC - {B^2}} \right|^{ - 2}}\left| {{{\left( { - A + AC - {B^2}} \right)}^3}} \right. \\
 \quad \,\,\,\,\,\,\,\,\,\,\,\,\,\left. { + {A^3}\left( {\left( { - A + AC - {B^2}} \right)\left( {1 + F} \right) + \left( {1 - C} \right){{\left( {1 + D} \right)}^2} + 2BE\left( {1 + D} \right) - A{E^2}} \right)} \right|. \\
 \end{array}
$$
By the fact $\varphi ' = {\alpha ^{ - 1}}{f^{ - 1}}{A^{{{ - 1}
\mathord{\left/
 {\vphantom {{ - 1} 2}} \right.
 \kern-\nulldelimiterspace} 2}}}$, we obtain
\begin{equation}\label{3.24}
\begin{array}{l}
 k_1^2 + k_2^2 = {\alpha ^{ - 2}}{f^{ - 2}}{A^{ - 3}}{\left| { - A + AC - {B^2}} \right|^{ - 2}}\left| {{{\left( { - A + AC - {B^2}} \right)}^3}} \right. \\
 \quad \quad \quad \;\left. {\,\, + {A^3}\left( {\left( { - A + AC - {B^2}} \right)\left( {1 + F} \right) + \left( {1 - C} \right){{\left( {1 + D} \right)}^2} + 2BE\left( {1 + D} \right) - A{E^2}} \right)} \right|. \\
 \end{array}
\end{equation}
and
\begin{equation}\label{3.25}
\begin{array}{l}
{k_1} = {\alpha ^{ - 1}}{f^{ - 1}}{A^{ - {3 \mathord{\left/
 {\vphantom {3 2}} \right.
 \kern-\nulldelimiterspace} 2}}}{\left( { - A + AC - {B^2}} \right)^{{1 \mathord{\left/
 {\vphantom {1 2}} \right.
 \kern-\nulldelimiterspace} 2}}}.
 \end{array}
\end{equation}
According to our assumption,
$$
\begin{array}{l}
 f = {\left( {1 + {g^2} + {h^2}} \right)^{{{ - 3} \mathord{\left/
 {\vphantom {{ - 3} 2}} \right.
 \kern-\nulldelimiterspace} 2}}}{\left| { - 1 - {g^2} - {h^2} + {{\dot g}^2} + {{\dot h}^2} + {{\left( {\dot gh - g\dot h} \right)}^2}} \right|^{{{ - 5} \mathord{\left/
 {\vphantom {{ - 5} 2}} \right.
 \kern-\nulldelimiterspace} 2}}} \\
 \,\,\,\,\,\,\,\,\,\,\,\left| {{{\left( { - 1 - {g^2} - {h^2} + {{\dot g}^2} + {{\dot h}^2} + {{\left( {\dot gh - g\dot h} \right)}^2}} \right)}^3}} \right. \\
 \,\,\,\,\,\,\,\,\,\,\,\,\,\, - {\left( {1 + {g^2} + {h^2}} \right)^3}\,\left( {{{\left( {g - \ddot g} \right)}^2} + {{\left( {h - \ddot h} \right)}^2}} \right. \\
 \quad \left. {\left. {\,\,\,\,\,\,\,\,\, - {{\left( {\left( {g\dot h - \dot gh} \right) + \left( {\dot g\ddot h - \ddot g\dot h} \right)} \right)}^2} + {{\left( {g\ddot h - \ddot gh} \right)}^2}} \right)} \right|, \\
 \end{array}
$$
we obtain
$$
\begin{array}{l}
 f = {A^{ - {3 \mathord{\left/
 {\vphantom {3 2}} \right.
 \kern-\nulldelimiterspace} 2}}}{\left| { - A + AC - {B^2}} \right|^{{{ - 5} \mathord{\left/
 {\vphantom {{ - 5} 2}} \right.
 \kern-\nulldelimiterspace} 2}}}\left| {{{\left( { - A + AC - {B^2}} \right)}^3}} \right. \\
 \,\,\,\,\,\,\,\,\left. {\, + {A^3}\left( {\left( {1 + F} \right) + \left( {1 - C} \right){{\left( {1 + D} \right)}^2} + 2BE\left( {1 + D} \right) - A{E^2}} \right)} \right|. \\
 \end{array}
$$
If we substitute the above equations (3.25) and (3.26), we obtain
$$
\begin{array}{l}
{k_1} = \alpha \left( {k_1^2 + k_2^2} \right).
 \end{array}
$$
The proof is completed.

\noindent In the above equation ${\mu _1} = {\mu _2} =  - 1$ which
is the special Case 1. This formula is the parametric equation of
generalized spacelike Mannheim curve with timelike second binormal
vector in the Minkowski space-time $E_1^4$.


\begin{thebibliography}{9}


\bibitem{Onei} B. O'Neill, {\it Semi-Riemannian Geometry with Applications to Relativity}, Academic Press, New York, (1983).
\bibitem{Str}  D. J. Struik, {\it Differential geometry}, Second ed., Addison-Wesley, Reading, Massachusetts, (1961).
\bibitem{Bal2} H. Balgetir, M. Bekta\d{s}, J. Inoguchi, {\it Null Bertrand curves in Minkowski 3-space and their characterizations}, Note Math., 23, no. 1, (2004).
\bibitem{Bal1} H. Balgetir, M. Bekta\d{s}, M. Erg\"{u}t, {\it Bertrand Curves for Nonnull Curves in 3-Dimensional Lorentzian Space}, Hadronic Journal, 27, (2004).
\bibitem{Liu}  H. Liu, F. Wang, {\it Mannheim Partner curves in 3-space}, Journal of Geometry, 88, 120-126, (2008).
\bibitem{Mat1} H. Matsuda, S. Yorozu, {\it Notes on Bertrand curves}, Yokohama Math. J. 50, no. 1-2, 41-58, (2003).
\bibitem{Mat2} H. Matsuda, S. Yorozu, {\it On generalized Mannheim curves in Euclidean 4-space}, (English), Nihonkai Math. J., 20, no. 1, 33-56, (2009).
\bibitem{Orb} K. Orbay, E. Kasap, {\it On Mannheim Partner Curves in ${E^3}$}, International Journal of Physical Sciences Vol. 4 (5), pp.
261-264, May, (2009).
\bibitem{Eis} L.P.A. Eisenhart, {\it Treatise on the Differential Geometry of Curves and Surfaces}, New York, Dover, (1960).
\bibitem{Car}  M.P. Do Carmo, {\it Differential Geometry of Curves and Surfaces}, Pearson Education, (1976).
\bibitem{Tig} O. Tigano,  {\it Sulla determinazione delle curve di Mannheim} , Matematiche Catania 3, 25-29, (1948).
\bibitem{Ekm} N. Ekmekci, K. Ilarslan, {\it On Bertrand curves and their characterization}, Differ. Geom. Dyn. Syst.(electronic), vol. 3, no. 2, (2001).
\bibitem{Blum} R. Blum, {\it A remarkable class of Mannheim curves}, Canad. Math. Bull. 9, 223-228, (1966).
\bibitem{Kuh}  W. Kuhnel, {\it Differential geometry, Curves-surfaces-manifolds}, Braunschweig, Wiesbaden, (1999).
\bibitem{Nad} Z. N\'{a}den\'{i}k,  {\it Bertrand curves in five-dimensional space}, (Russian), Czechoslovak Mathematical Journal, vol. 2, issue 1,  pp. 57-87, (1952).


\end{thebibliography}
\end {document}